\documentclass{amsart}
\usepackage[dvips]{graphicx}
\usepackage{amscd}
\usepackage{amsmath}
\usepackage{amsxtra}
\usepackage{amsfonts}
\usepackage{amssymb}
\newtheorem{theorem}{Theorem}[section]

\newtheorem{lemma}[theorem]{Lemma}

\theoremstyle{definition}
\newtheorem{definition}[theorem]{Definition}
\newtheorem{remark}[theorem]{Remark}

\theoremstyle{remark}

\renewcommand{\theclaim}{\textup{\theclaim}}

\newtheorem*{acknowledgements}{Acknowledgements}

\numberwithin{equation}{section}

\input cyracc.def

\begin{document}

\title[On the Laplacian subalgebra of certain tensors of group von Neumann algebras]{On the Laplacian subalgebra of certain tensors of group von Neumann algebras}
\author{Teodor \c Stefan B\^ ildea}
\address{Work done while at the Department of Mathematics\\
University of Iowa.\\
Currently at the Computational Biomedicine Lab, University of Houston} \email{teodor.bildea@mail.uh.edu}


\begin{abstract}
We consider the Laplacian (radial) von Neumann subalgebra of $\mathcal{L}(G)^{\overline{\otimes}_k}$, $k\ge 1$,
for two classes of groups: $G=\mathbb{F}_N$, the free group on $N\ge 2$ generators, and  $G=G_1\ast...\ast G_m$, the free product of $m\ge 3$ groups, each finite of order $p\ge 2$, $m\ge p$. In both cases the conditional expectation onto the Laplacian subalgebra is an asymptotic homomorphism, and therefore the Laplacian subalgebra is a  strongly singular masa.
\end{abstract}

\maketitle

\section{\label{Intro}Introduction} Let $M$ be  a type ${\rm II}_1$ factor with
trace $\tau$, and let $A\subset M$ be a maximal abelian von Neumann subalgebra
(masa) of $M$. Consider the group of normalizing unitaries of $A$,
$N_M(A)=\{ u\in A| u\mbox{ unitary },uAu^{\ast}=A\}$.  According to the size
of  $B=N_M(A)^{"}$ in $M$, two extreme situations can occur: $B=M$, that is $B$
has enough unitaries to generate $M$ and in this case $A$ is called {\bf regular}
or ${\bf Cartan}$;  $B=A$, in which case the  only normalizing unitaries are the
unitaries in $A$, and $B$ is called {\bf singular}. Cartan masas appear in a natural
way in cross-products. On the other hand, until recently, singular masas where hard
to find. R\u adulescu \cite{FR} showed, using the Pukansky invariant,  that the Laplacian
(or radial) subalgebra of $\mathcal{L}(\mathbb{F}_N)$, $N\ge 2,$ is a singular masa.
Together with Boca \cite{BR}, they showed that the same is true for the Laplacian
subalgebra of $\mathcal{L}(G)$, where $G=G_1\ast ...\ast G_m$, the free product of
$m\ge 3$ finite groups of order $p\ge2$, with $m\ge p$ .
\par The class of strongly singular  masas, introduced by Sinclair and Smith in  \cite{ss2}
(see also \cite{ss}, \cite{GSS}), turn out to be more tractable.
Using the notion of asymptotic homomorphism and good criteria to
detect such maps as well as strong singularity in general, examples
where produced demonstrating strong singularity for certan abelian
algebras. Among the first examples of strongly singular masas was
the Laplacian subalgebra of the free group factor.
\par The aim of this paper is to present  examples of strongly singular masas living inside tensor products of $\mathcal{L}(G)$ with itself, where $G$ is either a free group of rank at most $2$, or a free product of $m\ge3$ finite groups of order $p\ge 2$, with $m\ge p$. The examples are natural generalizations of the Laplacian subalgebra of $\mathcal{L}(\mathbb{F}_N)$, generated by the selfadjoint element
\[ w=\sum_{i=1}^n(g_i+g_i^{-1}),\]
where we have identified the generators $g_1,...,g_N$ of $\mathbb{F}_N$ with the unitaries in $\mathcal{L}(\mathbb{F}_N)$ given by the left regular representation. In $\mathcal{L}(\mathbb{F}_N)^{\overline{\otimes}_k}$ we consider  the abelian von Neumann subalgebra generated by
\[ w^{(k)}=\sum_{i=1}^n(g_i^{\otimes_k}+(g_i^{-1})^{\otimes_k}).\]
The investigation was started in \cite{bildea}, where we proved that this algebra is a strongly singular masa for each $k\ge 1$. We did this using the criteria introduced by Sinclair and Smith, which provides a new and easy way to prove that the conditional expectation onto this subalgebra is an asymptotic homomorphism. By the results of Sinclair and Smith, it follows that the respective subalgebra is a strongly singular masa. We will recall the results and point out similarities to the case of the free product of finite groups of the same order (but not necessarily isomorphic). We would like to thank Stuart White for sharing some of his dissertation results, essential for this paper.\\
The article is organized as follows: section \ref{pre} contains notation and some preliminaries. In section \ref{lemmas} we prove an essential lemma about counting words in a free product of finite groups. In section \ref{ssm} we review the notion of strongly singular masa and the tools needed to prove our result.  The last section contains the proof of the main result.

\section{\label{pre}Preliminaries and notation}

Let $K$ be any countable group.
The left regular representation is given by $\lambda:K\to B(l^2(K))$, $(\lambda(g)(\xi))(h)=\xi({g^{-1}h})$ for all $h\in K$, for all $\xi \in l^2(K)$ and for all $g\in K$.
The von Neumann algebra generated by $\{\lambda(g)|g\in K\}$ is denoted by $\mathcal{L}(K)$ and is called the von Neumann group algebra generated by $K$.
It is well known that if a countable group $K$ is an ICC group (all the conjugacy classes are infinite), then $\mathcal{L}(K)$ is a type ${\rm II}_1$ factor acting standardly on $l^2(K)$ and that with this identification the $\|\cdot\|_{\tau}$ norm coincides with the usual norm $\|\cdot\|_2$ on $l^2(K)$.
We will use the same notation $\|\cdot\|_2$ for the norm $\|\cdot\|_{\tau}$ given by $\|x\|_{\tau}=(\tau(x^{\ast}x)^{1/2})$.
We will identify elements in $K$ with their left regular representations. We will also identify the elements of $\mathcal{L}(K)$ with the corresponding vectors in $L^2(\mathcal{L}(K),\tau)=l^2(K)$.

Consider $G=G_1\ast ...\ast G_m$, the free product of $m\ge 3$ finite groups of order $p\ge2$, $m\ge p$. Let $G_i^{\circ}=G_i\setminus \{e_{G_i}\}$ for $i=1,2,...,m$. Each non-identity element can be uniquely written as $g=g_1...g_l$ where each $g_j$ is a member of some $G_i^{\circ}$ and no two consecutive $g_j$ lie in the same $G_i^{\circ}$. The length of $g$ is defined to be $l$ and this is denoted by $|g|=l$; the length of the empty word, or identity,  $e$, is set to be 0. This particular countable group $G$ is an ICC group, so $\mathcal{L}(G)$ is a type ${\rm II_1}$ factor.

\par If $H$ is a Hilbert space, we denote by $H^{\otimes_k}$ the $k$-folded tensor product of $H$ with itself (with the standard convention $H^{\otimes_0}:=\mathbb{C}$, $H^{\otimes_1}:=H$). An element $v\in H$ embeds as $v\otimes ...\otimes v=:v^{\otimes_k}\in H^{\otimes_k}$. For a von Neumann algebra $M$, $M^{\overline{\otimes}_k}$ denotes the $k$-folded von Neumann algebra tensor product of $M$ with itself ($M^{\overline{\otimes}_0}:=\mathbb{C}$, $M^{\overline{\otimes}_1}:=M$). For $x\in M$,  $x\otimes...\otimes x=:x^{\otimes_k}\in M^{\overline{\otimes}_k}$.

For $n\ge 0$ we define $w^{(k)}_n\in \mathcal{L}(G)^{\overline{\otimes}_k}$ to be the sum of all $k$-folded tensors of reduced words in $G$ of length $n$:
\[ w^{(k)}_0=e^{\otimes_k},\]
\[ w^{(k)}_n=\sum_{|v|=n}v^{\otimes_k}.\]
The following relations are known \cite{trenholme1},\cite{stuart} for $k=1$ and hold for any $k\ge 1$:
\begin{eqnarray}\label{wn}
(w_1^{(k)})^2=w_2^{(k)}-(p-2)w_1^{(k)}-m(p-1)w_0^{(k)};\\
 w^{(k)}_1w^{(k)}_n=w^{(k)}_nw^{(k)}_1=w^{(k)}_{n+1}+(p-2)w_n^{(k)}+(m-1)(p-1)w^{(k)}_{n-1},
\end{eqnarray}
for all $n\ge2$. These relations show that the von Neumann algebra $\mathcal{B}^{(k)}$ generated by $w^{(k)}_1$ contains each $w^{(k)}_n$, and hence is the weak closure of the span of these elements:
\[ \mathcal{B}^{(k)}=\overline{Sp\{{w^{(k)}_n|n\ge 0}\}}^{ w}\subset \mathcal{L}(G)^{\overline{\otimes}_k}.\]
We call $\mathcal{B}^{(k)}$ the {\bf Laplacian} or {\bf radial subalgebra}. In \cite{trenholme2} it is shown that for $k=1$ this algebra is a masa if and only if $m\ge p$. The selfadjoints $w^{(k)}_n$ are pairwise orthogonal with respect to the trace $\tau^{\otimes_k}$, and a simple counting argument shows that :
 \begin{equation}\label{norma wn}
  \|w^{(k)}_n\|^2_2=m(p-1)[(m-1)(p-1)]^{n-1},\quad n\ge 1,
  \end{equation}
  which represents the number of distinct words of length $n$.
  By the previous remarks
  \[ \{\frac{w^{(k)}_n}{\|w^{(k)}_n\|_2}|n\ge 0\} \]
   is an orthonormal basis for $L^2(\mathcal{B}^{(k)},\tau)$.
   Moreover, if $\mathbb{E}_k:\mathcal{L}(G)^{\overline{\otimes}_k}\to \mathcal{B}^{(k)}$ is the unique  trace-preserving conditional expectation, then
 \begin{equation}\label{exp x}
 \mathbb{E}_k(x)=\sum_{n=0}^{\infty}\tau^{\otimes_k}(xw^{(k)}_n)\frac{w^{(k)}_n}{\|w^{(k)}_n\|^2_2}.
 \end{equation}
 In particular, when $v$ is a reduced word of length $l$:
 \begin{equation}\label{vdekori}
  \mathbb{E}_k(v^{\otimes_k})=\frac{w^{(k)}_l}{\|w^{(k)}_l\|^2_2}.
  \end{equation}

\begin{lemma}\label{nonzero} For $k\ge 1$ let $x_1,...,x_{k}\in G$.
With the above notation ,
\[ \mathbb{E}_k(x_1\otimes...\otimes x_k)\ne 0\Leftrightarrow x_1=...=x_k.\]
\end{lemma}
{\it Proof: }
Recall from (\ref{exp x}) that
\[ \mathbb{E}_k(x_1\otimes...\otimes x_k)= \sum_{n=0}^{\infty}\tau^{\otimes_k}(x_1\otimes...\otimes x_kw^{(k)}_n)\frac{w^{(k)}_n}{\|w^{(k)}_n\|^2_2}.\]
Therefore
$\mathbb{E}_k(x_1\otimes...\otimes x_k)\ne 0$
if and only if there is some $n\ge 0$ so that :
\[\tau^{\otimes_k}(x_1\otimes...\otimes x_kw^{(k)}_n)\ne 0.\]
Because of the definitions of $w_n^{(k)}$ and $\tau^{\otimes_k}$ this can happen if and only if for some $n\ge 0$:
\[  \sum_{|v|=n}\tau^{\otimes_k}(x_1v\otimes...\otimes x_kv)\ne 0,\]
or
\[ \sum_{|v|=n}\tau(x_1v)...\tau(x_kv)\ne 0.\]
But this is true exactly when for some $v$ of length $n\ge 0$ we have
\[ x_1=...=x_k=v^{-1}. \quad \mbox{\rule{5pt}{5pt}}\]

\section{\label{lemmas}Counting words}
In this section we prove a technical lemma about solutions of a conjugacy equation in $G=G_1\ast ...\ast G_m$, the free product of $m\ge 2$ groups. A similar result holds for $\mathbb{F}_N$, the free group of rank $N\ge 2$ \cite{bildea}.
If $a,b\in G$ are two words such that their concatenation is a reduced word, we will denote this by $a\cdot b$. In general, the product (concatenation) will be denoted $ab$, meaning that cancellations may or may not take place.
\begin{lemma}\label{mylemma}
Let $a,b\in G$ be non-trivial elements.\\
(i) The equation $x\cdot a=b\cdot x$ has  in $G$ at most one solution of fixed positive length $l\ge 1$.
\\ (ii) The equation $xa=bx$ has  in $G$ at most one solution of fixed positive length $l>|a|+|b|$.
\end{lemma}
{\it Proof:} The proof is essentially the same as for the similar result that holds for $\mathbb{F}_N$, $N\ge2$; see \cite{bildea} for details.
\rule{5pt}{5pt}

\section{\label{ssm}Strongly singular MASA's}
A maximal abelian selfadjoint subalgebra $\mathcal{A}$ in a ${\rm II}_1$ factor $\mathcal{M}$ is singular if any unitary $u\in \mathcal{M}$ which normalizes $\mathcal{A}$ (i.e. $u\mathcal{A}u^{\ast}=\mathcal{A}$), must lie in $\mathcal{A}$.

To define the notion of strongly singular masa, consider a linear map $\phi:\mathcal{M}_1\to \mathcal{M}_2$ between two type ${\rm II}_1$ factors. There are several norms for $\phi$, depending on the norms considered on the two algebras. When $\mathcal{M}_1$ has the operator norm and $\mathcal{M}_2$ has the $\|\cdot\|_2$-norm given by the trace, we denote the resulting norm for $\phi$ by $\|\phi\|_{\infty,2}$, following \cite{ss2}.
\begin{definition} Suppose $\mathcal{A}$ is a masa  in a type ${\rm II}_1$ factor $\mathcal{M}$.\\
1)$\mathcal{A}$ is called $\alpha$-strongly singular (or simply strongly singular for $\alpha=1$) if
\[ \|\mathbb{E}_{u\mathcal{A}u^{\ast}}-\mathbb{E}_{\mathcal{A}}\|_{\infty,2}\ge \alpha\|u-\mathbb{E}_{\mathcal{A}}(u)\|_2,\]
for all unitaries $u\in \mathcal{M}$.\\
2) The conditional expectation $\mathbb{E}_{\mathcal{A}}$ is an asymptotic homomorphism if there is a unitary $u\in \mathcal{A}$ such that
\[ \lim_{|k|\to \infty}\|\mathbb{E}_{\mathcal{A}}(xu^ky)-\mathbb{E}_{\mathcal{A}}(x)\mathbb{E}_{\mathcal{A}}(y)u^k\|_2=0\]
for all $x,y\in \mathcal{M}$.
\end{definition}

  The following provides the criteria we will use in the last section to prove the main result of this paper.
\begin{theorem}\label{criteria}(Sinclair and Smith) Let $\mathcal{A}$ be a abelian von Neumann subalgebra of a type ${\rm II}_1$ factor $(\mathcal{M},tr)$, and suppose that there is a $\ast$-isomorphism $\pi:\mathcal{A}\to L^{\infty}[0,1]$ which induces an isometry from $L^2(\mathcal{A},tr)$ onto $L^2[0,1]$. Let $\{v_n\in \mathcal{A}|n\ge 0\}$ be an orthonormal basis for $L^2(\mathcal{A},tr)$, and let $Y\subseteq \mathcal{M}$ be a set whose linear span is norm dense in $L^2(\mathcal{M},tr)$. Let $\mathbb{E}_{\mathcal{A}}:\mathcal{M}\to \mathcal{A}$ be the unique conditional expectation satisfying $tr \circ \mathbb{E}_{\mathcal{A}}=tr$. If
\[ \sum^{\infty}_{n=0}\| \mathbb{E}_{\mathcal{A}}(xv_ny)-\mathbb{E}_{\mathcal{A}}(x)\mathbb{E}_{\mathcal{A}}(y)v_n\|_2^2<\infty\]
for all $x,y\in Y$, then $\mathbb{E}_{\mathcal{A}}$ is an asymptotic homomorphism and $\mathcal{A}$ is a strongly singular masa.
\end{theorem}

 \begin{remark}\label{r1}
  As abelian von Neumann algebras, $\mathcal{B}^{(k)}\subset B(L^2(\mathcal{B}^{(k)},\tau^{\otimes_k}))$ and $\mathcal{B}^{(1)}\subset B( L^2(\mathcal{B}^{(1)},\tau))$ are $\ast$-isomorphic and the isomorphism induces an isometry between the $L^2$-spaces, sending $w^{(k)}_n$ to $w^{(1)}_n, n\ge 0$. It is known \cite{trenholme1} that $\mathcal{B}^{(1)}$ is a masa, so it automatically satisfies the hypothesis of the above theorem. In \cite{stuart} it is shown that $\mathbb{E}_1:\mathcal{L}(B)\to \mathcal{B}^{(1)}$ satisfies the hypothesis of the above theorem (so $\mathbb{E}_1$  is an asymptotic homomorphism and $\mathcal{B}^{(1)}$ is a strongly singular masa). Via the $\ast$-isomorphism between $\mathcal{B}^{(k)}$ and $\mathcal{B}^{(1)}$ it follows that $\mathcal{B}^{(k)}$ also satisfies the hypothesis of the theorem, for all $k\ge 2$.
  \end{remark}
We further exploit the $\ast$-isomorphism between $\mathcal{B}^{(k)}$ and  $\mathcal{B}^{(1)}$ in the following :
\begin{lemma}\label{k0}Let $x,y\in G$ with $|x|=l_1\ge 2,|y|=l_2\ge 2$. The following equality holds for every $k\ge 1,n\ge l_1+l_2$:
\begin{equation}\label{eq1}
\|\mathbb{E}_k(x^{\otimes_k}w^{(k)}_ny^{\otimes_k})-\mathbb{E}_k(x^{\otimes_k})\mathbb{E}_k(y^{\otimes_k})w^{(k)}_n\|^2_2=\|\mathbb{E}_1(xw^{(1)}_ny)-\mathbb{E}_1(x)\mathbb{E}_1(y)w^{(1)}_n\|^2_2.
\end{equation}

\end{lemma}
 {\it Proof: } The proof involves the same computations as in the proof of the similar result that holds in the case of the free group of rank $N\ge 2$; see \cite{bildea} for details.\rule{5pt}{5pt}
\section{\label{fn}The case of tensors of $\mathcal{L}(G)$}

In what follows $G=G_1\ast ...\ast G_m$, the free product of $m\ge 3$ groups of order $p\ge 2$, $m\ge p$. Note that the groups $G_i, 1\le i\le m$ are not necessarily isomorphic. We will identify elements $g\in G$ with their left regular representations $\lambda(g)\in \mathcal{L}(G)$ and we will use the notation from section \ref{Intro}.

The next technical lemma contains the calculations needed to prove the main result. The proof is similar to the proof of the main result in \cite{bildea}.
\begin{lemma} Let $k\ge 1$ and $w^{(k)}_n=\sum_{|v|=n}v^{\otimes_k}$. Let $x_1,...,x_{k},y_1,...,y_{k}$ be words in $G$. Then:
\begin{equation}\label{ineq}
\sum_{n=0}^{\infty}\frac{1}{\|w^{(k)}_n\|^2_2}\| \mathbb{E}_k(x_1\otimes...\otimes x_{k}w^{(k)}_ny_1\otimes...\otimes y_{k})-
\end{equation}
\[ -\mathbb{E}_k(x_1\otimes...\otimes x_{k})\mathbb{E}_k(y_1\otimes...\otimes y_{k})w^{(k)}_n\|_2^2<\infty.
\]
\end{lemma}

{\it Proof: } We look first at the case $\mathbb{E}_k(x_1\otimes...\otimes x_k)\mathbb{E}_k(y_1\otimes...\otimes y_k)\ne 0$. By lemma \ref{nonzero}, it follows that $x_1=...=x_k=:x$ and $y_1=...=y_k=:y$. Let $l_1=|x|,l_2=|y|$. If $l_1=1$ or $l_2=1$ the inequality (\ref{ineq}) is trivial. Indeed, suppose $l_1=1$ (the other case is even easier). Then we actually have
\[ \|\mathbb{E}_k(x^{\otimes_k}w^{(k)}_ny^{\otimes_k})-\mathbb{E}_k(x^{\otimes_k})\mathbb{E}_k(y^{\otimes_k})w^{(k)}_n\|^2_2=\]
\[\| \mathbb{E}_k(w^{(k)}_ny^{\otimes_k})-\mathbb{E}_k(y^{\otimes_k})w^{(k)}_n\|^2_2=\]
\[ \|w^{(k)}_n\mathbb{E}_k(y^{\otimes_k})-\mathbb{E}_k(y^{\otimes_k})w^{(k)}_n\|^2_2=0,\]
because $\mathcal{B}^{(k)}$ is abelian. We are left with the case $l_1,l_2\ge 2$. Using lemma \ref{k0}, the finiteness of the sum in (\ref{ineq}) follows in this case from the corresponding result in \cite{ss}.

For remainder of the proof assume  \[\mathbb{E}_k(x_1\otimes...\otimes x_k)\mathbb{E}_k(y_1\otimes...\otimes y_k)= 0.\] By the lemma \ref{nonzero}, there must be at least two distinct $x_{i_1},x_{i_2}$ and at least two distinct $y_{j_1},y_{j_2}$. The nonzero terms in the series (\ref{ineq})
\[\sum_{n=0}^{\infty}\frac{1}{\|w^{(k)}_n\|^2_2}\| \mathbb{E}_k(x_1\otimes...\otimes x_kw^{(k)}_ny_1\otimes...\otimes y_k)\|_2^2=\]

 \[\sum_{n=0}^{\infty}\frac{1}{\|w^{(k)}_n\|^2_2}\| \sum_{|v|=n}\mathbb{E}_k(x_1vy_1\otimes...\otimes x_kvy_k)\|_2^2\]
are given by the solutions $v$ of length $n\ge 0$  of the equations:
\[ x_1vy_1=x_2vy_2=...=x_kvy_k.\]
Observe that a solution exists only if $x_i=x_j\Leftrightarrow y_i=y_j$ for $1\le i,j\le k+1$. Without loss of generality, assume $x_1\ne x_2(\Leftrightarrow y_1\ne y_2)$.
Any solution of length $n$ for the above system of equations is a solution of $v(y_1y^{-1}_2)=(x^{-1}_1x_2)v$. Clearly $x^{-1}_1x_2\ne e\ne  y_1y^{-1}_2$. By lemma \ref{mylemma}, for each fixed $n> {\rm max}\{|x_i| : 1\le i\le k+1\}+{\rm max}\{|y_j|:1\le j\le k+1\}=:n_0$ there is at most one solution to this equation. So
\[\sum_{n=n_0}^{\infty}\frac{1}{\|w^{(k)}_n\|^2_2}\| \sum_{|v|=n}\mathbb{E}_k(x_1vy_1\otimes...\otimes x_kvy_k)\|_2^2\le\]
\[\sum_{n=n_0}^{\infty}\frac{1}{\|w^{(k)}_n\|^2_2}\| \frac{w^{(k)}_n}{\|w^{(k)}_n\|^2_2}\|^2_2=\]
\[\sum_{n=n_0}^{\infty}\frac{1}{\|w^{(k)}_n\|^4_2} \stackrel{(\ref{norma wn})}{<}\infty.\qquad\rule{5pt}{5pt}\]

\begin{theorem}
For each $k\ge 1$,the conditional expectation $\mathbb{E}_k$ onto the Laplacian subalgebra $\mathcal{B}^{(k)}$ of $\mathcal{L}(G)^{\overline{\otimes}_k}$, generated by
\[ w^{(k)}_1=\sum_{v\in G,|v|=1}v^{\otimes_k}, \]
is an asymptotic homomorphism and $\mathcal{B}^{(k)}$ is a strongly singular masa.
\end{theorem}
{\it Proof: } For $k=1$, $\mathcal{L}(G)^{\overline{\otimes}_1}=\mathcal{L}(G)$ and the result was proved in \cite{stuart}. For $k\ge 2$, apply  theorem \ref{criteria}: use remark \ref{r1} and the proof for $k=1$ to verify the statement about the existence of the $\ast$-isometry; use $Y=\{x_1\otimes ...\otimes x_{k}| x_i\in G,1\le i\le k\}$ and $v_n:=w^{(k)}_n/\|w^{(k)}_n\|_2,n\ge 0$. The theorem is a direct consequence of the above technical lemma and theorem \ref{criteria}.
\rule{5pt}{5pt}
\begin{acknowledgements} . I want to first thank the organizers of GPOTS 2005 for inviting me and for giving me the opportunity to present my results. I also  want to thank Stuart White for valuable comments and for sharing some results from his Ph.D. thesis.
\end{acknowledgements}


\begin{thebibliography}{20}
\bibitem{bildea} B\^ ildea, T.\c S. {\it The Laplacian subalgebra of $\mathcal{L}(\mathbb{F}_N)^{\overline{\otimes}_k}$ is a strongly singular MASA}, preprint
\bibitem{BR} Boca, F.; R\u adulescu, F. {\it Singularity of radial subalgebras in ${\rm II}_1$ factors asociated with free products of groups}, J. Funct. Anal., {\bf 103}  (1992), no. 1, 138--159.
\bibitem{trenholme1}Cohen, J.M; Trenholme, A.R. {\it Orthogonal polynomials with a constant recursion formula and an application to harmonic analysis}, J. Funct. Anal.,  {\bf 59} (1984), pp 175-184.
\bibitem{FR}R\u adulescu, F. {\it Singularity of the radial subalgebra of $\mathcal{L}(\mathbb{F}_N)$ and the Puk\'ansky invariant},  Pacific J. Math. , {\bf 151}, no.2 (1991), pp 297-306.
\bibitem{GSS}Robertson, G.; Sinclair, A. M.; Smith, R. R. {\it Strong singularity for subalgebras of finite factors}, Intern. J. Math., {\bf 14}, No. 3 (2003) 235--258.
\bibitem{ss} Sinclair, A. M.; Smith, R. R. {\it The Laplacian MASA in a free group factor} Trans. Amer. Math. Soc. {\bf 355} (2003), no. 2, 465--475.
\bibitem{ss2} Sinclair, A. M.; Smith, R. R. {\it Strongly singular masas in type ${\rm II}_1$ factors} Geom. Funct. Anal. {\bf 12} (2002), no. 1, 199--216.
\bibitem{trenholme2}Trenholme, A.R. {\it Maximal abelian subalgebras of function algebras associated with free products}, J. Funct. Anal., {\bf 79} (1988), pp 342-350.
\bibitem{stuart} White, S. {\it Strong Singularity of Radial Subalgebras in ${\rm II}_1$ factors associated with free products of finite groups} - private communication
\end{thebibliography}
 \end{document}